\documentclass[a4paper, 10pt, conference]{ieeeconf}      

\IEEEoverridecommandlockouts                              

\usepackage[a4paper, left=1.91cm, right=1.31cm, top=1.91cm,bottom=3.67cm]{geometry}

\usepackage{cite}
\usepackage{amsmath,amssymb,amsfonts}
\usepackage{algorithmic}
\usepackage{graphicx}
\usepackage{textcomp}
\usepackage{xcolor}

\usepackage{paralist}
\usepackage{hyperref}
\usepackage{comment}

\newcommand{\rk}{\mathrm{rk}\,}

\newcommand{\norm}[1]{\left\Vert #1 \right\Vert}

\newcommand{\JK}[1]{\textcolor{blue}{#1}}

\newtheorem{theorem}{Theorem}[section]

\newtheorem{example}[theorem]{Example}
\newtheorem{definition}[theorem]{Definition}
\newtheorem{lemma}[theorem]{Lemma}
\newtheorem{definition-proposition}[theorem]{Definition-Proposition}
\newtheorem{remark}[theorem]{Remark}

\newtheorem{algorithm}[theorem]{Algorithm}

\def\BibTeX{{\rm B\kern-.05em{\sc i\kern-.025em b}\kern-.08em
    T\kern-.1667em\lower.7ex\hbox{E}\kern-.125emX}}
\begin{document}

\title{Funnel MPC for port-Hamiltonian systems\\
\thanks{This work was funded by the Deutsche Forschungsgemeinschaft (DFG, German Research Foundation) - Project-ID 544702565.}
}

\author{Loris K\"onitzer$^{1}$ and Jonas Kirchhoff$^{1}$ and Thomas Berger$^{1}$
\thanks{$^{1}$L. K\"onitzer, J. Kirchhoff and T. Berger are with the Institute of Mathematics, Martin Luther University Halle-Wittenberg, Halle (Saale), Germany
        {\tt\small \{loris-benedikt.koenitzer, jonas.kirchhoff, thomas.berger\}@mathematik.uni-halle.de}}%
}

\maketitle
\thispagestyle{empty}
\pagestyle{empty}

\begin{abstract}
Port-Hamiltonian systems are an energy-based modelling formalism, that has strong ties to the physical models and allows interconnection of systems. Funnel MPC, a Model Predictive Control (MPC) scheme, allows for output tracking with prescribed error bounds. Earlier works focused on systems with stable internal dynamics, which typically excludes port-Hamiltonian systems. We introduce a modified stage cost design that takes into account the port-Hamiltonian structure of the system and allows deactivation of the funnel cost. We show initial and recursive feasibility for a class of port-Hamiltonian systems with relative degree one including unstable internal dynamics.
\end{abstract}
\section{Introduction}
 
The port-Hamiltonian systems formalism is a well-established framework for modelling and control of physical systems. For a comprehensive overview, we refer the reader to~\cite{SchaftJeltsema14,JacoZwar14}.

One appealing feature of port-Hamiltonian systems is that the formalism is intrinsically based on power-flow between the (virtual ideal) components of the system; e.g., a real-existing spring is virtually decomposed into an ideal spring, modelled as an energy storage component which stores and releases potential energy, and a damper, where energy is dissipated due to internal friction of the spring. In particular, the energy supplied to the system appears naturally as product of the external port variables. This property has prompted~\cite{FaulMascPhilSchaWort21,FaulMascPhilSchaWort22} to propose this quantity as the natural cost function for port-Hamiltonian optimal control problems.

Here we are interested in the problem of reference tracking with transient safety guarantees. The safe region is given by a time-dependent Euclidean ball around the reference trajectory, and the control objective becomes the design of a controller that guarantees that the trajectory remains, for all time, within the safe region. Globally, the problem is solved by the funnel controller, cf.~\cite{BergIlchRyan25}. Recently, the idea of funnel control was combined with model predictive control (MPC) giving the funnel MPC scheme~\cite{BergDennIlchWort22,BergDenn24}, which produces reference tracking with transient guarantees.

As MPC is an optimal control formalism, this prompts the combined control objective: Derive an MPC scheme for port-Hamiltonian systems which achieves reference tracking with transient safety guarantees and minimal supplied energy.

Port-Hamiltonian systems are in some sense a difficult system class for reference tracking, as it is not uncommon that the internal dynamics are unstable, which is in large parts owed to the physics of the system.

In this note, we propose a combined MPC scheme for the relative degree one linear-quadratic prototype of port-Hamiltonian systems. We show reference tracking for systems with stable zero dynamics, which is no restriction for relative degree one port-Hamiltonian systems, and discuss the consequences for the Funnel MPC algorithm. We shall explore a novel stage cost design, that incorporates the port-Hamiltonian structure and allows the deactivation of funnel costs, when far from its boundaries. We show initial and recursive feasibility and demonstrate that the control objective is satisfied.

\section{Reminder on port-Hamiltonian systems}

We consider the linear port-Hamiltonian system

\begin{equation}
\begin{aligned}\label{eq:lq_pH}
    \tfrac{\mathrm{d}}{\mathrm{d}t}x & = (J-R)Q x+Bu,\\
    y & = B^\top Qx,
\end{aligned}
\end{equation}
where $J = -J^\top\in\mathbb R^{n\times n}$, $R = R^\top\in\mathbb R^{n\times n}$ is non-negative, $Q=Q^\top$ is non-negative and $B\in\mathbb{R}^{n\times m}$. A solution of~\eqref{eq:lq_pH} is $(x,y,u)\in\mathcal{C}(I,\mathbb{R}^n)\times L^2_{\mathrm{loc}}(I,\mathbb{R}^{2m})$, $I\subseteq\mathbb{R}$ an open interval, so that, for all $t_0,t\in I$,
\begin{align*}
    x(t) &= x(t_0)+\int_{t_0}^{t}(J-R)Q x(\tau)+Bu(\tau)\,\mathrm{d}\tau
\end{align*}
and $y = B^\top Qx$. The pair $(u,y)$ are the \textit{port variables} of~\eqref{eq:lq_pH}. The Hamiltonian function $H(x):=\frac{1}{2}x^\top Qx$ models the energy of the system.

\begin{example}
    Consider a force-controlled mass-spring-damper system. Suppose a mass $M$ is attached to a spring (assumed to have negligible mass) that may only move in horizontal direction and has a spring constant $k>0$. Let $q$ denote the position of the mass and $p$ its momentum. The energy of the system is given by the sum of the potential and kinetic energy:
    \begin{equation*}
        H(q,p)=\frac{k}{2}q^2+ \frac{1}{2M}p^2.
    \end{equation*}
    With a force input, the port-Hamiltonian system is given by
    \begin{equation}\label{eq:massspringdamper}
    \begin{aligned}
        \frac{\mathrm{d}}{\mathrm{d}t}\begin{pmatrix}
            q\\p
        \end{pmatrix}&=\begin{bmatrix}
            0&\frac{1}{M}\\
            -k&-\frac{r}{M}
        \end{bmatrix}\begin{pmatrix}
            q\\p
        \end{pmatrix}
        +\begin{bmatrix}
            0\\1
        \end{bmatrix}u,\\
        y&=\frac{1}{M}p,
    \end{aligned}        
    \end{equation}
    where $r\geq0$ is a damping coefficient and $y$ is the collocated velocity output. That means the product of the input and output has the physical interpretation as the supplied power of the system.
\end{example}

A straightforward calculation gives the well-known energy balance equation for the port-Hamiltonian system.

\begin{lemma}
Let $I\subseteq\mathbb{R}$ be an open interval, $t_0< t_1\in I$ and $(x,y,u)\in\mathcal{C}(I,\mathbb{R}^n)\times L^2_{\mathrm{loc}}(I,\mathbb{R}^{2m})$ be a solution of~\eqref{eq:lq_pH}. Then,
\begin{equation}\label{eq:energy_balance}
\begin{aligned}
H(x(t_1)) & = H(x(t_0)) - \int_{t_0}^{t_1}\norm{R^{\frac{1}{2}}Qx(t)}^2\,\mathrm{d}t\\
&\qquad + \int_{t_0}^{t_1}u(t)^\top y(t)\,\mathrm{d}t,
\end{aligned}
\end{equation}
where $\norm{\cdot}$ denotes the Euclidean norm, and $R^{\frac{1}{2}}$ is the matrix square root of $R$.
\end{lemma}

The change in energy along the trajectories of~\eqref{eq:lq_pH} is therefore governed by the internal dissipation of energy $\int_{t_0}^{t_1}\norm{R^{\frac{1}{2}}Qx}^2\,\mathrm{d}t$ and the supplied energy $\int_{t_0}^{t_1}u^\top y\,\mathrm{d}t$ as the integral of the supplied power.

Finally, for simplicity of exposition, we restrict ourselves to the case that~\eqref{eq:lq_pH} has \textit{relative degree one}. 

\begin{definition}[cf.~\cite{HinrPritIII}]\label{def:rel_deg}
    The system (\ref{eq:lq_pH}) has relative degree one if, and only if, $B^\top QB$ is invertible.
\end{definition}

It is not assumed that the Hamiltonian function is strictly convex. For reference tracking, the input-to-output behaviour is most relevant. We demonstrate that the input-to-output behaviour of~\eqref{eq:lq_pH} is the input-to-output behaviour of a port-Hamiltonian system with strictly convex Hamiltonian function.

\begin{lemma}\label{lem:Qisinvertible}
There exists a coordinate transformation $Sx=(x_1,x_2)\in\mathbb{R}^{(n-l)+l}$, where $\rk Q=l$, and matrices $A\in\mathbb{R}^{(n-\ell)\times\ell}$, $\widehat{J},\widehat{R},\widehat{Q}\in\mathbb{R}^{\ell\times\ell}$ with $\widehat{J}+\widehat{J}^\top = \widehat{R}-\widehat{R}^\top = 0$, $\widehat{Q}>0$ and $\widehat{R}\geq 0$, $\widehat{B}\in\mathbb{R}^{\ell\times m}$ and $C\in\mathbb R^{(n-l)\times m}$ which converts~\eqref{eq:lq_pH} into
\begin{align}
\tfrac{\mathrm{d}}{\mathrm{d}t}x_1 & = Ax_2+Cu\label{eq:singulardynamics}\\
\begin{split}
\tfrac{\mathrm{d}}{\mathrm{d}t}x_2 & = (\widehat{J}-\widehat{R})\widehat Qx_2+ \widehat{B}u\\
y & = \widehat{B}^\top \widehat{Q}x_2
\end{split}\label{eq:nonsingulardynamics}
\end{align}
\eqref{eq:lq_pH} has relative degree one if, and only if, $\rk\widehat{B} = m$.
\end{lemma}
\begin{proof}
    There exists an orthogonal matrix $S\in\mathbb R^{n\times n }$, such that $ S QS^\top=\mathrm{diag}(0,\widehat Q)$, where $\widehat Q\in\mathbb R^{l\times l}$ and $\widehat Q>0$ and hence
    \begin{equation*}
        S(J-R)QS^\top=S(J-R)S^\top S Q S^\top= S(J-R)S^\top\begin{bmatrix}
            0&0\\0&\widehat{Q}
        \end{bmatrix}.
    \end{equation*}
    Decompose accordingly
    \begin{align*}
        S JS^\top&=\begin{bmatrix}
            J_{1,1}&J_{1,2}\\
            J_{2,1}&J_{2,2}
        \end{bmatrix}\in\mathbb R^{(n-l)+l\times (n-l)+l},\\
        S RS^\top&=\begin{bmatrix}
            R_{1,1}&R_{1,2}\\
            R_{2,1}&R_{2,2}
        \end{bmatrix}\in\mathbb R^{(n-l)+l\times (n-l)+l}
    \end{align*}
    and put $\widehat{J}:=J_{2,2}$, $\widehat{R}:=R_{2,2}$ and $A:=(J_{1,2}-R_{1,2})\widehat Q$. Then $\widehat{J} = -\widehat{J}^\top$ and $\widehat{R} = \widehat{R}^\top\geq 0$. Likewise, define
    \begin{equation*}
        SB=\begin{bmatrix}
            C\\
            \widehat{B}
        \end{bmatrix}\in\mathbb{R}^{(n-l)+l\times m}.
    \end{equation*}
    By Definition~\ref{def:rel_deg}, $B^\top QB$ is invertible. Invoking
    \begin{equation*}
        B^\top Q B=[C^\top~\widehat{B}^\top]\begin{bmatrix}
            0&0\\0&\widehat{Q}
        \end{bmatrix}\begin{bmatrix}
            C\\
            \widehat{B}
        \end{bmatrix}=\widehat{B}^\top \widehat{Q}\widehat{B},
    \end{equation*}
    and invertibility of $\widehat{Q}$ gives $\rk\widehat{B} = m$ if, and only if,~\eqref{eq:lq_pH} has relative degree one. Lastly:
    \begin{align*}
        y&=B^\top Qx=B^\top S^\top SQS^\top\begin{pmatrix}
            x_1\\x_2
        \end{pmatrix}\\
        &=[
            C^\top~ \widehat{B}^\top]\begin{bmatrix}
            0&0\\0&\widehat{Q}
        \end{bmatrix}\begin{pmatrix}
            x_1\\x_2
        \end{pmatrix}=\widehat{B}^\top\widehat{Q}x_2
    \end{align*}
\end{proof}

We have demonstrated that the state-space of~\eqref{eq:lq_pH} can be split into~\eqref{eq:singulardynamics} and~\eqref{eq:nonsingulardynamics} so that the input-output behaviours of~\eqref{eq:lq_pH} and~\eqref{eq:nonsingulardynamics} coincide and~\eqref{eq:nonsingulardynamics} is port-Hamiltonian with strictly convex Hamiltonian function. We may therefore restrict the analysis to~\eqref{eq:nonsingulardynamics} and assume w.l.o.g. that~\eqref{eq:lq_pH} and~\eqref{eq:nonsingulardynamics} coincide.

A second decomposition of systems, which is routinely used for proofs in funnel control and funnel MPC, is the so-called \textit{Byrnes-Isidori} normal form, cf. \cite{ByrnesIsidori}. For the port-Hamiltonian system~\eqref{eq:nonsingulardynamics}, the Byrnes-Isidori normal form is structure-preserving in the following way.

\begin{lemma}\label{lem:structure_preserving_BI}
Consider~\eqref{eq:lq_pH} with $Q>0$ and $\rk B = m$. Let the columns of $W$ form a basis of $\ker B^\top$. The matrix
$\Phi:=\left[\begin{smallmatrix}
B^\top Q\\W^\top
\end{smallmatrix}\right]$
transforms~\eqref{eq:lq_pH} \JK{\eqref{eq:nonsingulardynamics}} into the port-Hamiltonian system
\begin{equation}\label{eq:structure_preserving_BI}
\begin{aligned}
\tfrac{\mathrm{d}}{\mathrm{d}t}\begin{pmatrix}
y\\\eta
\end{pmatrix} = \Phi(J-R)\Phi^\top\begin{pmatrix}
\Gamma^{-1}y\\
(W^\top Q^{-1}W)^{-1}\eta
\end{pmatrix} + \begin{bmatrix}
\Gamma\\0
\end{bmatrix}u,
\end{aligned}
\end{equation}
where $\Gamma := B^\top QB$ is the high-frequency gain matrix. In particular, the zero-dynamics is the port-Hamiltonian system
\begin{equation}\label{eq:zero_dynamics}
\tfrac{\mathrm{d}}{\mathrm{d}t}\eta = W^\top(J-R)W(W^\top Q^{-1}W)^{-1}\eta.
\end{equation}
\end{lemma}
\begin{proof}
By definition of $W$, $W^\top B = 0$, and plugging in $\Phi$ gives $\Phi B = \left[\begin{smallmatrix}
B^\top QB\\W^\top B
\end{smallmatrix}\right] = \left[\begin{smallmatrix}
\Gamma\\0
\end{smallmatrix}\right]$. Further,
\begin{align*}
\Phi Q^{-1}\Phi^\top & = \begin{bmatrix}
B^\top QB & B^{\top}W \\W^\top Q^{-1}QB & W^\top Q^{-1}W
\end{bmatrix}\\ & = \begin{bmatrix}
B^\top QB & 0\\
0 & W^\top Q^{-1}W
\end{bmatrix}.
\end{align*}
 In particular $\Phi Q^{-1}\Phi^\top$ is invertible. Moreover $\Phi$ has an inverse given by:
 \begin{equation*}
     \Phi^{-1}= Q^{-1}\Phi^\top \left(\Phi Q^{-1}\Phi^\top\right)^{-1}.
  \end{equation*}
With $\Phi^{-\top} Q\Phi^{-1} = (\Phi Q^{-1}\Phi^\top)^{-1}$,~\eqref{eq:lq_pH} becomes \eqref{eq:structure_preserving_BI}. 
\end{proof}

We write
\begin{align*}
    L&=W^\top (J-R)\Phi^\top \left(\begin{smallmatrix}
   0\\ (W^\top Q^{-1}W)^{-1}
\end{smallmatrix}\right),\\
K&=W^\top (J-R)\Phi^\top \left(\begin{smallmatrix}
   \Gamma^{-1}\\0
\end{smallmatrix}\right).
\end{align*}
The internal dynamics of (\ref{eq:lq_pH}) are then given by
\begin{equation}\label{eq:internaldynamics}
    \dot \eta(t)=L\eta(t)+Ky(t).
\end{equation}
 In funnel MPC and funnel control it is typically assumed that the internal dynamics of (\ref{eq:lq_pH}) are stable (i.e. the $\eta$ is stable in the sense that it is linearly bounded by the initial value and the supremum norm of the output $y$). As the following example will illustrate, for many port-Hamiltonian systems this requirement is oftentimes not fulfilled.

\begin{example}
If we choose a unit mass $M=1$ the mass-spring-damper system~\eqref{eq:massspringdamper} is in Byrnes-Isidori form~\eqref{eq:structure_preserving_BI}. The internal dynamics is the elongation $q$, which grows to infinity for any constant positive momentum. The internal dynamics are unstable.
\end{example}
Note that the damping does not affect the conclusion that the system has unstable internal dynamics, as damping is typically a force. There are of course port-Hamiltonian systems that have stable non-trivial internal dynamics, but what this example reinforces is that for a general port-Hamiltonian system it is neither guaranteed nor a trivial assumption for the internal dynamics to be stable. Instead we will rely on the zero dynamics of the system.

\begin{lemma}\label{lem:quasistable}
The zero-dynamics~\eqref{eq:zero_dynamics} is Lyapunov-stable.
\end{lemma}
\begin{proof}
By~\eqref{eq:energy_balance}, $\eta\mapsto\frac{1}{2}\big(\eta^\top(W^\top Q^{-1}W)^{-1}\eta\big)$ is a Lyapunov function for~\eqref{eq:zero_dynamics}.
\end{proof}
 
\section{Problem formulation}

Our objective is for the output $y$ to track a reference trajectory $y_{\mathrm{ref}}$ within a prescribed time-varying error bound $\psi>0$. To be precise, for all $t\geq 0$ we want to ensure that
\begin{align*}
    \Vert e(t) \Vert<\psi(t),\quad 
    e:=y-y_{\mathrm{ref}},
\end{align*}
where $\psi\in W^{1,\infty}([0,\infty),\mathbb{R})\cap L^1([0,\infty),\mathbb{R})$, $\psi>0$ and $y_{\mathrm{ref}}\in W^{1,\infty}([0,\infty),\mathbb{R}^m)$. We also impose the following integrability condition on $y_{\mathrm{ref}}$: there exists some $D>0$ such that for all $t\geq 0$
\begin{equation}\label{eq:yrefbound}
    \left\Vert \int _0^t e^{L(t-\tau)}Ky_{\mathrm{ref}}(\tau)\mathrm{d}\tau\right\Vert\leq D.
\end{equation}
\begin{remark}
    A typical assumption on the funnel would be a positive lower bound leading only to practical tracking, cf. \cite{BergDennIlchWort22}. We explicitly exclude this when we assume the funnel to be in $L^1$. We will show that we find uniform bound on the control with no lower bound on the funnel.\\ 
    The condition (\ref{eq:yrefbound}) also deviates from the standard funnel MPC literature. It is a special case of the trackability assumption, that was explored in \cite[(A2)]{BERGERnonmimimumphase} in the context of funnel control. For the mass-spring-damper system (\ref{eq:massspringdamper}) we can calculate that $L=0$ and $K=1$. The integral (\ref{eq:yrefbound}) then simplifies to:
    \begin{equation*}
        \left|\int _0^t y_{\mathrm{ref}}(\tau)\mathrm{d}\tau\right|\leq D,
    \end{equation*}
    essentially saying the position curve traced out by the velocity reference should stay bounded. This condition is not surprising: Had we chosen a position output of the system (which would not have been port-Hamiltonian), the system has relative degree two and stable internal dynamics. Boundedness of the position reference signal (which is what this integral represents) is a standard assumption for relative degree two funnel MPC. So the bound (\ref{eq:yrefbound}) excludes cases where a constant non-zero velocity is tracked, but allows for references such as those that tend to zero sufficiently fast (allowing stabilization) or are periodic (with mean at 0). With that being said, when dealing with funnel control and robust funnel MPC as introduced in \cite{RobustFMPC} the consequence of choosing a funnel boundary with no positive lower bound will lead to unbounded control values, nullifying whatever advantage we might have gained. 
\end{remark}

\section{Funnel MPC and energy supply stage costs}
The typical stage cost function $\tilde l_\psi:\mathbb R_{\geq 0}\times\mathbb R^m\times\mathbb R^m\rightarrow\mathbb R\cup\{\infty\}$ of a funnel MPC algorithm is given by 
\begin{align}\label{eq:FMPC_cost_function}
    \tilde l_\psi(t,y,u)&=\frac{\Vert y-y_{\mathrm{ref}}(t)\Vert^2}{\psi(t)^2-\Vert y-y_{\mathrm{ref}}(t)\Vert^2}+\lambda\Vert u\Vert^2
\end{align}
when $\Vert y-y_{\mathrm{ref}}(t)\Vert<\psi(t)$, and $\tilde l_\psi(t,y,u)=\infty$ otherwise, where $\lambda\geq 0$ is a design parameter.
This choice of stage cost achieves that $y$ stays in the funnel by penalizing close proximity to the funnel boundary (similar to a barrier function from optimization). The quadratic term in $u$ has the effect that we can minimize the control while remaining in the funnel. 
In the port-Hamiltonian framework the supplied energy is a quantity naturally arising from the dissipation inequality (\ref{eq:energy_balance}). It is then also natural to try to minimize the supplied energy in the context of optimal control as seen in \cite{FaulMascPhilSchaWort21},\cite{FaulKirMehrPhilSchaWort23}.
This prompted us to propose a novel stage cost that seeks to modify the typical funnel stage cost with a term that minimizes the supplied energy. We also modify the funnel costs such that they only activate close to the funnel boundary. For some $a\in[0,1)$ let an activation function $\alpha:\mathbb R_{\geq0}\times\mathbb R^m\to\mathbb{R}$ be defined by
\begin{equation*}
    \alpha: (t,e) \mapsto \begin{cases}
        1,&\Vert e\Vert>a\psi(t),\\
        0,&\Vert e\Vert\leq a\psi(t).
    \end{cases}
\end{equation*}
If we set $a=0$ the funnel is always active. While not actively discussed in this note, we may also use a continuous activation function to ease potential numerical burdens. We define $l_\psi:\mathbb R_{\geq 0}\times\mathbb R^m\times\mathbb R^m\rightarrow\mathbb R\cup\{\infty\}$ with 
\begin{equation}\label{eq:pHFMPC_cost_function}
\begin{aligned}
    l_\psi(t,y,u)&=\alpha(t,e)\frac{\Vert e\Vert^2}{\psi(t)^2-\Vert e\Vert^2}+\lambda u^\top y,\ e = y-y_{\mathrm{ref}}(t)
\end{aligned}
\end{equation}
when $\Vert y-y_{\mathrm{ref}}(t)\Vert<\psi(t)$ and $l_\psi(t,y,u)=\infty$ otherwise.
\begin{remark}
    We note that all theoretic results presented in this work, that apply to the stage cost $l_\psi$, also apply to the typical $\tilde l_\psi$ stage cost. These results are interesting in their own right, as showing the feasibility funnel MPC for system with unstable internal dynamics has not been shown so far. 
\end{remark}

\begin{algorithm}[FMPC]\label{alg:FMPCalgorithm}~\\
\textbf{Given:} System \eqref{eq:lq_pH}, reference signal $y_{\mathrm{ref}},$ initial state $x_0\in\mathbb R^n$ with error $\Vert e(0)\Vert<\psi(0)$, control bound $\Theta>0$ and the stage cost $l_\psi$.\\
\textbf{Set:} the time shift $\delta>0$, the prediction horizon $T\geq \delta$ and initial time $\hat{t}:=0$.\\
\textbf{Steps:}
\begin{enumerate}
    \item Obtain a measurement of the state at $\hat t$ and set $\hat x=x(\hat{t} )$.
    \item Compute a solution $u^*\in L^\infty([\hat t,\hat t+T],\mathbb R^m)$ of the optimal control problem
    \begin{equation}\label{eq:OCP}
    \begin{aligned}
            &\min_{\Vert u\Vert_\infty\leq \Theta}\int_{\hat t}^{\hat t+T}l_\psi(t,y(t),u(t))\mathrm{d}t\\
            \mathrm{s.t.}&\ (x,y,u) ~\mbox{is a solution of (\ref{eq:lq_pH}) and}~x(\hat t)=\hat x.
    \end{aligned}
\end{equation}
    \item Apply the feedback law
    \begin{equation}\label{eq:Feedbacklaw}
        \mu:[\hat t,\hat t+\delta)\times \mathbb R^n\rightarrow\mathbb R^m,~\mu(t,\hat x)=u^*(t)
    \end{equation}
    to system \eqref{eq:lq_pH}. Increase $\hat t$ by $\delta$ and go to Step 1).
\end{enumerate}
\end{algorithm}
\section{Initial and recursive feasibility}
    To show that the algorithm \ref{alg:FMPCalgorithm} is correct we need to show initial and recursive feasibility, which in practice means we need to show that the optimal control problem \eqref{eq:OCP} has a solution at each step. To that end, for $\hat t\ge 0$, consider the cost functional 
    \begin{equation}\label{eq:costfunctional}
        \begin{aligned}
            J_T^\psi(\cdot,\hat t,\cdot): L^\infty([\hat t,\hat t+T],\mathbb R^m)\times\mathbb R^n\rightarrow\mathbb R\cup\{\infty\}\\
(u,x_0)\mapsto\int_{\hat t}^{\hat t +T}l_\psi(t,y(t),u(t))\mathrm{d}t,
        \end{aligned}
        \end{equation}
where 
    \begin{equation*}
     y (t) \!=\!B^\top Q\left(e^{(J-R)Q(t-\hat t)}x_0\!+\!\int_{\hat t}^{t}e^{(J-R)Q(t-\tau)}Bu(\tau)\mathrm{d}\tau\right).
        \end{equation*}
    The following lemma verifies, that if a feasible control exists, then there exist a optimal control.
\begin{lemma}\label{lem:OCP}
    Consider the optimal control problem (\ref{eq:OCP}) with control bound $\Theta>0$ and prediction horizon $T>0$. Let $\hat t\ge 0$, $x_0\in\mathbb{R}^n$ and suppose there exists some $u\in  L^\infty([\hat t,\hat t+T],\mathbb R^m)$ such that $J_T^\psi(u,\hat t,x_0)<\infty$ and $\Vert u\Vert_\infty\leq \Theta$. Then there exists some $u^*\in  L^\infty([\hat t,\hat t+T],\mathbb R^m)$ with $\Vert u^*\Vert_\infty\leq \Theta$ that minimizes~(\ref{eq:OCP}), meaning
    \begin{equation*}
        J_T^\psi(u^*,\hat t,x_0)=\inf_{\Vert u\Vert_\infty\leq \Theta} J_T^\psi(u,\hat t,x_0).
    \end{equation*}
\end{lemma}
\begin{proof}
    Without loss of generality we may assume that $\hat t=0$ and we write $J_T^\psi(u,x_0):=J_T^\psi(u,0,x_0)$. By the energy balance (\ref{eq:energy_balance}) the cost functional is bounded from below:
    \begin{align*}
        J_T^\psi(u,x_0)&=\int_{0}^{T}l_\psi(t,y(t),u(t))\mathrm{d}t\\
        &=\int_{0}^{T}\underbrace{\alpha(t,e(t))\frac{\Vert e(t)\Vert^2}{\psi(t)^2-\Vert e(t)\Vert^2}}_{\geq0}+\lambda u(t)^\top y(t)\,\mathrm{d}t\\
        &\geq \lambda(H(x(T))-H(x(0)))\\
        &\geq -\lambda H(x_0).
    \end{align*}
    Hence, there exists a minimizing sequence $(u_n)_{n\in\mathbb N}\subset L^\infty([0,T],\mathbb R^m)$  such that 
    \begin{align*}
        -\infty<\inf_{\Vert u\Vert_\infty\leq \Theta} J_T^\psi(u,x_0)
        =\lim_{n\rightarrow\infty} J_T^\psi(u_n,x_0)<\infty.
    \end{align*}
    Since $[0,T]$ is compact, $L^\infty([0,T],\mathbb R^m)\subset  L^2([0,T],\mathbb R^m)$, $u_n$ is a bounded sequence in the Hilbert space $L^2([0,T],\mathbb R^m)$. Hence, there exist a weakly convergent subsequence (which we will not relabel) $u_n\rightharpoonup u^*$ in $L^2([0,T],\mathbb R^m)$. The set of all controls $\{u\in L^\infty([0,T],\mathbb R^m)|\Vert u\Vert_\infty\leq \Theta\} $ is closed and convex and therefore weakly sequentially closed \cite{OCodPDEs}, so that $\Vert u^*\Vert_\infty\leq \Theta$.
    The associated sequence of state trajectories $(x_n)_{n_\in\mathbb N}\subset \mathcal{C}([0,T],\mathbb R^n)$ is defined by
    $$x_n(t)=e^{(J-R)Qt}x_0+\int_0^t e^{(J-R)Q(t-\tau)}Bu_n(\tau)\mathrm{d}\tau$$
    for all $t\in[0,T]$ and $n\in\mathbb N$. This sequence is uniformly bounded, since for all $n\in\mathbb N$
    \begin{equation*}
    \Vert x_n(t)\Vert\leq\Vert e^{(J-R)Qt}\Vert \Vert x_0\Vert+T\sup_{t\in[0,T]}\Vert e^{(J-R)Qt}\Vert \Vert B\Vert \Theta,
    \end{equation*}
    and uniformly equicontinuous since for all $n\in\mathbb N$ and $t_1,t_2\in[0,T]$,
    \begin{align*}
    \Vert x_n(t_1)-x_n(t_2)\Vert\leq \int_{t_1}^{t_2}\Vert (J-R)Q x_n(t)+Bu_n(t)\Vert\mathrm{d}t\\
    \leq \left(\Vert (J-R)Q\Vert \sup_{t\in[0,T]}\Vert x_n(t)\Vert+\Vert B\Vert \Theta\right)|t_2-t_1|.
    \end{align*}
    By the Arzelà–Ascoli theorem there exists a subsequence of $(x_n)_{n\in\mathbb N}$ (which we will not relabel) that uniformly converges to some $x^*\in \mathcal{C}([0,T],\mathbb R^n)$. To show that this $x^*$ is the trajectory associated with the control $u^*$ first note that by the bounded convergence theorem
\begin{equation*}
    \int_0^t (J-R)Qx_n(\tau)\mathrm{d}\tau\rightarrow \int_0^t (J-R)Qx^*(\tau)\mathrm{d}\tau
\end{equation*}
    for $n\rightarrow\infty$. Now due to weak convergence of $u$ we have
\begin{equation*}
     \int_0^t Bu_n(\tau)\mathrm{d}\tau\rightarrow \int_0^t Bu^*(\tau)\mathrm{d}\tau
\end{equation*}
for $n\rightarrow\infty$ (we can choose the rows of $B$ as test functions).Therefore, $x^*$ is the solution to 
\begin{equation*}
    x^*(t)=x_0+\int_0^t(J-R)Qx^*(\tau)+Bu^*(\tau)\mathrm{d}\tau.
\end{equation*}
Likewise, define $y^*=B^\top Qx^*$. To show that $(u^*,y^*)$ is an optimal solution, we show that cost functional (\ref{eq:costfunctional}) is to be semicontinuous from below. First, we calculate (for brevity we omit the dependence on $t$)
\begin{align*}
    \left|\int_0^Tu^{*\top} y^*-u_n^\top y_n\mathrm{d}t\right|&\leq \int_0^T|u_n^\top( y^*-y_n)|\mathrm{d}t\\
    &+\left|\int_0^T(u^*-u_n)^\top y^*\mathrm{d}t\right|\rightarrow0
\end{align*}
for $n\rightarrow\infty$ due to boundedness of $u_n$, weak convergence of $u_n$ and uniform convergence of $y_n$.\\
When dealing with the funnel term let $t\in[0,T]$. Let 
\begin{equation*}
    e_n=y_n-y_{\mathrm{ref}},\quad e^*=y^*-y_{\mathrm{ref}}
\end{equation*}
be the errors associated with $y_n$ and $y^*$, respectively. Furthermore, we write
\begin{equation*}
    \gamma(t,e)=\begin{cases}
        \alpha(t,e)\frac{\Vert e\Vert^2}{\psi(t)^2-\Vert e\Vert^2},& \Vert e\Vert<\psi(t)\\
        \infty,&\Vert e\Vert\geq\psi(t)
    \end{cases}
\end{equation*}
for all $e\in\mathbb R^m$. For $\Vert y^*(t)-y_{\rm ref}(t)\Vert\leq a\psi(t)$:
\begin{align*}
    \gamma(t,e_n(t))\geq0=\gamma(t,e^*(t)).
\end{align*}
For $\psi(t)>\Vert y^*(t)-y_{\mathrm{ref}}(t)\Vert> a\psi(t)$ there exists a neighborhood around $e^*(t)$, where $\gamma$ is continuous and therefore
\begin{equation*}
    \gamma(t,e_n(t))\rightarrow \gamma(t,e^*(t)),\quad n\to\infty.
\end{equation*}
If $\Vert y^*(t)-y_{\rm ref}(t)\Vert\geq\psi(t)$, then $\gamma(t,e^*(t))=\infty$, By continuity $\gamma(t,e_n(t))$ diverges as well. Therefore
\begin{equation*}
    \liminf_{n\rightarrow\infty}\gamma(t,e_n(t))\geq \gamma(t,e^*(t)).
\end{equation*}
Applying Fatou's Lemma yields
\begin{align*}
    \liminf_{n\rightarrow\infty}\int _0^T\gamma(t,e_n(t))\mathrm{dt}&\geq\int _0^T\liminf_{n\rightarrow\infty}\gamma(t,e_n(t))\mathrm{dt}\\
    &\geq \int_0^T\gamma(t,e^*(t))\mathrm{d}t.
\end{align*}
Finally, recall that $u_n$ is minimizing sequence:
\begin{align*}
    \liminf_{n\rightarrow\infty} J_T^\psi(u_n,x_0)&=\liminf_{n\rightarrow\infty} \int_{0}^{T}l_\psi(t,y_n(t),u_n(t))\mathrm{d}t\\
    &=\liminf_{n\rightarrow\infty} \int_{0}^{T}\gamma(t,e_n(t))+\lambda u_n(t)^\top y_n(t)\mathrm{d}t\\
    &\geq \int_{0}^{T}\gamma(t,e^*(t))+\lambda u^{*\top}(t) y^*(t)\mathrm{d}t\\
    &= J_T^\psi(u^*,x_0)\\
    &\geq  \liminf_{n\rightarrow\infty} J_T^\psi(u_n,x_0),
\end{align*}
which verifies $ J_T^\psi(u^*,x_0)= \liminf_{n\rightarrow\infty} J_T^\psi(u_n,x_0)$.
\end{proof}
The following lemma establishes the equivalence between trajectories that stay within the funnel and trajectories with finite cost with respect to the cost functional (\ref{eq:costfunctional}).
\begin{lemma}\label{lem:equivsets}
    Let $(x,y,u)$ be a solution of (\ref{eq:lq_pH}) with $\Vert u\Vert_\infty<\infty$ and let $t_0\leq t_1$. Then the following statements are equivalent:\\
    (i) $\int_{t_0}^{t_1} l_\psi(t,y(t),u(t))\mathrm{d}t<\infty$.\\
    (ii) $\Vert e(t)\Vert<\psi(t)$ for all $t\in[t_0,t_1]$.
    \end{lemma}
\begin{proof}
    Since $\psi>0$, it has a positive lower bound on $[t_0,t_1]$. Also, $x$ is bounded on $[t_0,t_1]$. Moreover the activation function $\alpha$ can be neglected, since the singularity appears when $\Vert e\Vert>a\psi$. Under these considerations the assertion follows analogously to \cite[Thm. 4.3]{BergDennIlchWort22}, where we replace non-negativity with the boundedness from below shown in Lemma \ref{lem:OCP}. 
\end{proof}
The following lemma is essential when constructing the control bound $\Theta$. We show that for a fixed initial value the state trajectory remains bounded. 
\begin{lemma}\label{lem:boundedstate}
    Suppose $Q$ is invertible and $\mathrm{rk}B=m$. Let $(x,y,u)$ be a solution of (\ref{eq:lq_pH}) where $x(0)=x_0$ and $\Vert  y(t)-y_{\mathrm{ref}}(t)\Vert<\psi(t)$ for all $t\geq 0$. Then there exists $\theta>0$ such that for all $t\geq0$ we have $\Vert x(t)\Vert\leq \theta$.
\end{lemma}
\begin{proof}
    Without loss of generality we assume that the system is in Byrnes-Isidori form (\ref{eq:structure_preserving_BI}). The following inequalities hold:
    \begin{align*}
        \Vert y\Vert_\infty&\leq \Vert y_{\mathrm{ref}}\Vert_\infty+\Vert y-y_{\mathrm{ref}}\Vert_\infty\leq\Vert y_{\mathrm{ref}}\Vert_\infty+\Vert \psi\Vert_\infty.
    \end{align*}
    Consider the internal dynamics (\ref{eq:internaldynamics}), which for $\eta(0)=\eta_0$ have the solution
    \begin{equation*}
        \eta(t)=e^{Lt}\eta_0+\int_0^te^{L(t-\tau)}Ky(\tau)\mathrm{d}\tau.
    \end{equation*}
    From Lemma \ref{lem:quasistable} it follows that $t\mapsto e^{Lt}$ is bounded. Together with the integrability assumptions (\ref{eq:yrefbound}) we obtain:
    \begin{align*}
        \Vert\eta(t)\Vert &\leq\left\Vert e^{Lt}\eta_0\right\Vert +\left\Vert \int_0^te^{L(t-\tau)}Ky_{\mathrm{ref}}(\tau)\mathrm{d}\tau\right\Vert \\
        &+\left\Vert \int_0^te^{L(t-\tau)}Ke(\tau)\mathrm{d}\tau\right\Vert\\
        &\leq\left\Vert e^{L\,\cdot}\right\Vert_\infty\Vert\eta_0\Vert +D+\left\Vert e^{L\,\cdot}\right\Vert_\infty\Vert K\Vert \underbrace{\Vert e\Vert_{L^1}}_{<\Vert \psi\Vert_{L^1}}.
    \end{align*}
    Then choose $\theta>0$ such that:
    \begin{equation*}
        \Vert x(t) \Vert\leq \Vert \Phi^{-1}\Vert(\Vert y_{\mathrm{ref}}\Vert_\infty+\Vert \psi\Vert_\infty+\Vert\eta\Vert_\infty)\leq \theta 
    \end{equation*}
\end{proof}
We show Initial and recursive feasibility of the FMPC Algorithm \ref{alg:FMPCalgorithm}.
\begin{theorem}\label{Thm:InitialandRecursiveFeasibility}
    Consider system \eqref{eq:lq_pH} with a relative degree of one and initial value $x(0)=x_0$. Let  $y_{\mathrm{ref}}\in W^{1,\infty}([0,\infty),\mathbb{R}^m)$ fulfill condition (\ref{eq:yrefbound}) and $\psi\in W^{1,\infty}([0,\infty),\mathbb{R})\cap L^1([0,\infty),\mathbb{R})$, $\psi>0$. Let $\Vert B^\top Qx_0-y_{\mathrm{ref}}(0)\Vert<\psi(0)$. Then there exists $\Theta>0$ only dependent on $x_0. \psi$ and $y_{\mathrm{ref}}$ such that the FMPC Algorithm \ref{alg:FMPCalgorithm} with time shift $\delta>0$, prediction horizon $T\geq\delta$, and stage cost $l_\psi$
    is initially and recursively feasible, i.e., for $\hat t=0$ and for every successor step $\hat t\in \delta\mathbb N$ the  optimal control problem (\ref{eq:OCP}) has a solution. In particular, the closed-loop system consisting of (\ref{eq:lq_pH}) and the feedback law (\ref{eq:Feedbacklaw}) has a (not necessarily unique) global solution $(x,y,u)$, where $u$ is given by:
    \begin{equation*}
        u(t)=\mu(t,x(\hat t)),\quad t\in [\hat t,\hat t+\delta),\ \hat t\in\delta\mathbb N_0.
    \end{equation*}
    Furthermore this solution satisfies:\\
    (i) $\Vert u(t)\Vert\leq \Theta$ for all $t\geq 0$.\\
    (ii) The error $e=y-y_{\mathrm{ref}}$ fulfills $\Vert e(t)\Vert <\psi(t)$ for all $t\geq0$.
\end{theorem}
\begin{proof}
    As was shown in Lemma \ref{lem:Qisinvertible}, we may assume that $Q$ is invertible, since the singular dynamics of $Q$ have no influence on the Funnel MPC Algorithm \ref{alg:FMPCalgorithm}.
    Set $\Theta=\Vert\Gamma^{-1}\Vert\left(\left\Vert B^\top Q(J-R)Q\right\Vert \theta+\Vert \dot\psi\Vert_\infty+\Vert\dot y_{\mathrm{ref}}\Vert_\infty\right)$, where $\theta$ is the constant given by Lemma~\ref{lem:boundedstate} and $\Gamma$ is the high frequency gain matrix defined in Lemma \ref{lem:structure_preserving_BI}. Suppose the Funnel MPC Algorithm~(\ref{alg:FMPCalgorithm}) was successful on its first $N$ time steps, such that we are now at time $t_0=\delta N$, $N\in\mathbb{N}_0$. Let $(x,y,u)$ be a solution obtained by applying the feedback law~(\ref{eq:Feedbacklaw}) on $[0,t_0]$. For all $t\in[0,t_0]:$ $\Vert e(t)\Vert <\psi(t)$. We construct a feasible control to the OCP (\ref{eq:OCP}). Apply the state feedback
    \begin{equation*}
        u(t)=\Gamma^{-1}\left(-B^\top Q(J-R)Qx(t)+\frac{e(t_0)}{\psi(t_0)}\dot\psi(t)+\dot y_{\mathrm{ref}}(t)\right).
    \end{equation*}
    We can then calculate for $t\geq t_0$:
    \begin{align*}
       & \Vert e(t)\Vert=\left\Vert \int_{t_0}^t\dot y(\tau)-\dot y_{\mathrm{ref}}(\tau)\mathrm{d}\tau+e(t_0)\right\Vert\\
        &=\left\Vert \int_{t_0}^tB^\top Q(J-R)Qx(\tau)+\Gamma u(\tau)-\dot y_{\mathrm{ref}}(\tau)\mathrm{d}\tau+e(t_0)\right\Vert\\
         &=\left\Vert \int_{t_0}^t \frac{e(t_0)}{\psi(t_0)}\dot\psi(\tau)\mathrm{d}\tau+e(t_0)\right\Vert\\
         &=\left\Vert \frac{e(t_0)}{\psi(t_0)}(\psi(t)-\psi(t_0))+e(t_0)\right\Vert\\
         &=\frac{\Vert e(t_0)\Vert }{\psi(t_0)}\psi(t)<\psi(t).
    \end{align*}
    In other words, $e(t)$ stays inside the funnel for all $t\geq t_0$. Due to Lemma~\ref{lem:boundedstate}, $\Vert x(t)\Vert\leq \theta$ and therefore $\Vert u(t)\Vert\leq \Theta$ for all $t\geq t_0$. By Lemma \ref{lem:equivsets}, $u$ is then a feasible control and we can apply Lemma~\ref{lem:OCP} to show that there exists an optimal control $u^*$ for the OCP (\ref{eq:OCP}). The associated optimal output trajectory must, by Lemma~\ref{lem:equivsets}, again stay in the funnel, thus completing the proof.
\end{proof}

\begin{remark}
    By Lemma \ref{lem:boundedstate}, if $Q>0$, then the state $x$ remains bounded under application of the FMPC Algorithm~\ref{alg:FMPCalgorithm}.
\end{remark}

\begin{figure}[ht]
    \centering
    \includegraphics[width=\columnwidth]{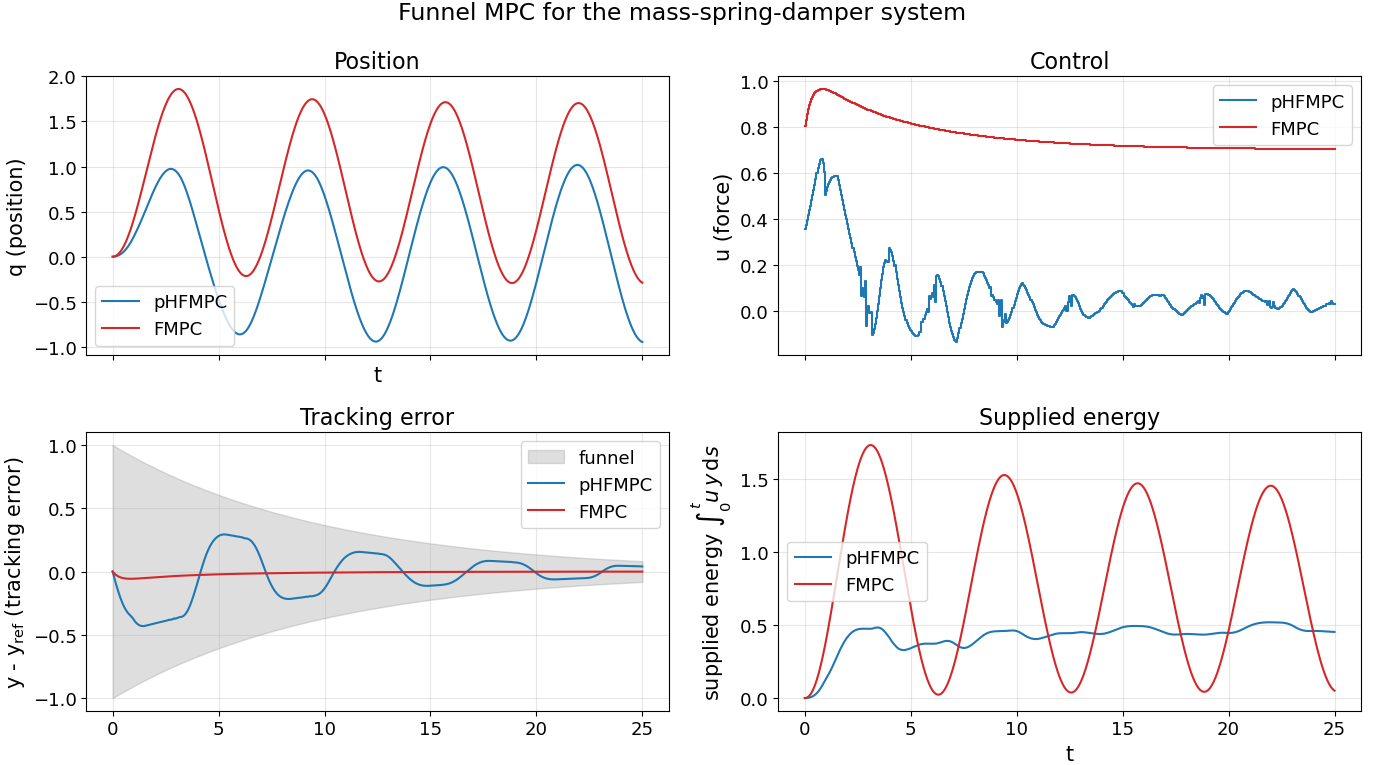}
    \caption{Simulation result with the reference tracking on the bottom-left and the supplied energy on the bottom-right.}
    \label{fig:simulation_result}
\end{figure}

\section{Numerical experiment}

In this section, we consider the mass-spring-damper system~\eqref{eq:massspringdamper}, and we compare the outcome of the MPC Algorithm \ref{alg:FMPCalgorithm} with the two stage costs $l_\psi$ (pHFMPC) and $\tilde l_\psi$ (FMPC). With unit mass, unit spring constant and vanishing damping, we obtain
\begin{equation}\label{eq:numerical_experiment}
\begin{aligned}
\tfrac{\mathrm{d}}{\mathrm{d}t}\begin{pmatrix}
q\\p
\end{pmatrix} & = \begin{bmatrix}
0 & 1\\-1 & 0
\end{bmatrix}\begin{pmatrix}
q\\p
\end{pmatrix} + \begin{pmatrix}
0\\1
\end{pmatrix}u\\
y & = p
\end{aligned}
\end{equation}
We consider a funnel radius $\psi(t)= e^{-0.1t}$, the reference $y_{\mathrm{ref}}(t) = \sin(t)$, and zero initial value. We set $\lambda=0.1$, $a=0.5$, $T=1$, $\delta=0.05$ and $\Theta=100$. The result of the experiment is presented in Figure~\ref{fig:simulation_result}. The resulting position trajectory can be seen on the top left. We can see that pHFMPC does not follow the reference as closely as the classic FMPC algorithm (bottom left). The benefits are reduced control values (top right) and smaller peaks in the supplied energy (bottom right). This can be explained by the pHFMPC trying to minimize the terminal stored energy of the system. Note that these effects are in part due to the activation function.

\section{Conclusion}

The port-Hamiltonian funnel MPC stage cost~\eqref{eq:pHFMPC_cost_function} combines the control objectives of minimal energy supply and of funnel control. It was observed that port-Hamiltonian systems are, in general, not treated by classical FMPC since the internal dynamics are not stable. Nonetheless, it was demonstrated that the FMPC algorithm is well-posed. A numerical example illustrates the different behaviours of classical and port-Hamiltonian FMPC.

\bibliographystyle{IEEEtran}
\bibliography{source}
\end{document}